\documentstyle[amscd,amssymb,verbatim,epsf]{amsart}

\begin{document}
\theoremstyle{plain}
\newtheorem{Thm}{Theorem}
\newtheorem{Cor}{Corollary}
\newtheorem{Con}{Conjecture}
\newtheorem{Main}{Main Theorem}
\newtheorem{Lem}{Lemma}
\newtheorem{Prop}{Proposition}

\theoremstyle{definition}
\newtheorem{Def}{Definition}
\newtheorem{Note}{Note}

\theoremstyle{remark}
\newtheorem{notation}{Notation}
\renewcommand{\thenotation}{}

\errorcontextlines=0
\numberwithin{equation}{section}
\renewcommand{\rm}{\normalshape}%

\title[Weingarten surfaces in hyperbolic 3-space]%
   {A characterization of Weingarten surfaces in hyperbolic 3-space}
\author{Nikos Georgiou}
\address{Nikos Georgiou\\
          Department of Computing and Mathematics\\
          Institute of Technology, Tralee\\
          Clash\
          Tralee\\
          Co. Kerry\\
          Ireland.}
\email{nikos.georgiou@@research.ittralee.ie}
\author{Brendan Guilfoyle}
\address{Brendan Guilfoyle\\
          Department of Computing and Mathematics\\
          Institute of Technology, Tralee\\
          Clash\\
          Tralee\\
          Co. Kerry\\
          Ireland.}
\email{brendan.guilfoyle@@ittralee.ie}

\date{8th March 2009}

\begin{abstract}
We study 2-dimensional submanifolds of the space ${\mathbb{L}}({\mathbb{H}}^3)$ of oriented geodesics of hyperbolic 
3-space, endowed with the canonical neutral K\"ahler structure. Such a surface 
is Lagrangian iff there exists a surface in ${\mathbb{H}}^3$ orthogonal to the geodesics of $\Sigma$. 

We prove that 
the induced metric on a Lagrangian surface in ${\mathbb{L}}({\mathbb{H}}^3)$  has zero Gauss curvature 
iff the orthogonal surfaces in ${\mathbb{H}}^3$ are Weingarten: the 
eigenvalues of the second fundamental form are functionally related. We then classify the totally null surfaces in
${\mathbb{L}}({\mathbb{H}}^3)$  and recover the well-known holomorphic constructions of flat and CMC 1  
surfaces in ${\mathbb{H}}^3$.
\end{abstract}

\maketitle

Recently the existence and uniqueness of a canonical neutral K\"ahler structure on the space 
${\mathbb{L}}({\mathbb{H}}^3)$ of oriented geodesics of hyperbolic 3-space ${\mathbb{H}}^3$  has been established
\cite{gag} \cite{Salvai}. The main purpose of this paper is to apply this 
work to the study of surfaces $S$ in ${\mathbb{H}}^3$. The oriented geodesics normal to $S$ form a surface in 
${\mathbb{L}}({\mathbb{H}}^3)$ which is Lagrangian with respect to this K\"ahler structure. In fact, the study of
surfaces in  ${\mathbb{H}}^3$ is equivalent, at least locally, to the study of Lagrangian surfaces in
${\mathbb{L}}({\mathbb{H}}^3)$.

Special classes of surfaces in hyperbolic 3-space have been studied for many decades, with various constructions being 
developed and applied to the particular class of surfaces under consideration. For example, surfaces of constant 
mean curvature 1
have been investigated in \cite{bryant}, while flat surfaces have been treated in \cite{umehara1}  \cite{umehara2}  
\cite{roitman}.
The common feature of these surfaces is that they are Weingarten \cite{weingarten}: they have some specified 
functional 
relationship between the eigenvalues of the second fundamental form. Other forms of this relationship have been considered 
\cite{eps} \cite{KaS} and uniqueness results obtained for Weingarten surfaces in ${\mathbb R}^3$ \cite{chern} \cite{haw}.

In this paper we give a new characterization of the Weingarten condition for surfaces in ${\mathbb{H}}^3$:

\vspace{0.2in}

\noindent {\bf Main Theorem}:

{\it
Let $S\subset{\mathbb{H}}^3$ be a $C^4$ smooth immersed oriented surface and $\Sigma\subset{\mathbb{L}}({\mathbb{H}}^3)$ 
be the Lagrangian 
surface formed by the oriented geodesics normal to $S$. Assume that the metric ${\mathbb G}_\Sigma$
induced on $\Sigma$ by the neutral K\"ahler metric is non-degenerate. 
Then $S$ is Weingarten iff the Gauss curvature of ${\mathbb G}_\Sigma$ is zero.
}

\vspace{0.2in}

The proof of this result follows from a careful study of 2-dimensional submanifolds of ${\mathbb{L}}({\mathbb{H}}^3)$
using local coordinates, moving frames and the correspondence space. In the following section we give the 
geometric background (further details can be found in \cite{gag}) and explore the submanifold theory of surfaces in 
${\mathbb{L}}({\mathbb{H}}^3)$. Section three contains the proof of the Main Theorem.

In the final section we classify the totally null surfaces in
${\mathbb{L}}({\mathbb{H}}^3)$  and recover the holomorphic constructions of flat and CMC 1  
surfaces in ${\mathbb{H}}^3$.

\vspace{0.2in}

\section{Geometric Background}

We briefly recall the basic construction of the canonical neutral K\"ahler metric on the space 
${\mathbb{L}}({\mathbb{H}}^3)$ of oriented geodesics of ${\mathbb{H}}^3$ - further details can be found in \cite{gag}.
We use one of two models of ${\mathbb H}^3$, the Poincar\'e ball model:
\[
 B^3=\{(y^{1},y^{2},y^{3})\in {\mathbb{R}}^3\;|\; (y^{1})^2+(y^{2})^2+(y^{3})^2< 1\},
\]
with standard coordinates $(y^{1},y^{2},y^{3})$ on ${\mathbb{R}}^3$, and hyperbolic metric
\[
ds^2=\frac{4[(dy^1)^{2}+(dy^2)^{2}+(dy^3)^{2}]}{[1-(y^{1})^2-(y^{2})^2-(y^{3})^2]^2}.
\]
and the upper-half space model:
\[
{\mathbb{R}}^3_{+}=\{\;(x^{0},x^{1},x^{2})\in {\mathbb{R}}^3 \;|\; x^{0}>0\;\},
\]
for standard coordinates $(x^{0},x^{1},x^{2})$ on ${\mathbb{R}}^3$. In these coordinates the hyperbolic metric 
has expression:
\[
d\tilde{s}^2=\frac{(dx^{0})^2+(dx^{1})^2+(dx^{2})^2}{(x^{0})^2}.
\]
These are related by the mapping 
${\mathbb{R}}^3_{+}\rightarrow  B^3\colon (x^{0},x^{1},x^{2}) \mapsto (y^{1},y^{2},y^{3})$
defined by
\[
y^{1}=\frac{2x^1}{(x^0+1)^2+(x^1)^2+(x^2)^2},\quad 
y^{2}=\frac{2x^2}{(x^0+1)^2+(x^1)^2+(x^2)^2},
\]
\[
y^{3}=\frac{(x^0)^2+(x^1)^2+(x^2)^2-1}{(x^0+1)^2+(x^1)^2+(x^2)^2}.
\]

An oriented geodesic in ${\mathbb{H}}^3$ is uniquely determined by its beginning and end point on the boundary of the 
ball model, and so ${\mathbb{L}}({\mathbb{H}}^3)$ can be identified with $S^2\times S^2-\Delta$, where $\Delta$ is the 
diagonal in $S^2\times S^2$. Endowing $S^2\times S^2-\Delta$ with the standard differentiable structure, a tangent 
vector to an oriented geodesic $\gamma\in {\mathbb{L}}({\mathbb{H}}^3)$ can then be identified with an orthogonal 
Jacobi field along $\gamma\subset {\mathbb{H}}^3$. 

Rotation of Jacobi fields through 90$^{0}$ about $\gamma$ defines an almost complex structure on 
${\mathbb{L}}({\mathbb{H}}^3)$. This almost complex structure is integrable, and so ${\mathbb{L}}({\mathbb{H}}^3)$ 
becomes a complex surface, which turns out to be biholomorphic to 
${\mathbb{P}}^1\times {\mathbb{P}}^1-\overline{\Delta}$. Here $\overline{\Delta}$ is the ``reflected'' diagonal: 
in terms of holomorphic coordinates $(\mu_1,\mu_2)$ on 
${\mathbb{P}}^1\times {\mathbb{P}}^1$, $\overline{\Delta}=\{(\mu_1,\mu_2):\mu_1\bar{\mu}_2=-1\}$. This 
distinction between ${\mathbb{P}}^1\times {\mathbb{P}}^1-\Delta$ and 
${\mathbb{P}}^1\times {\mathbb{P}}^1-\overline{\Delta}$ is crucial, as explained in section \ref{s:hol}.

The complex structure ${\mathbb{J}}$ on ${\mathbb{L}}({\mathbb{H}}^3)$ can be supplemented with a compatible symplectic 
structure $\Omega$, which has the following expression in holomorphic coordinates:
\begin{equation}\label{e:sympl}
\Omega=-\left[\frac{1}{(1+\mu_{1}\bar{\mu}_{2})^2}d\mu_{1}\wedge d\bar{\mu}_{2}
           +\frac{1}{(1+\bar{\mu}_{1}\mu_{2})^2}d\bar{\mu}_{1}\wedge d\mu_{2}\right].
\end{equation}
Together we obtain a K\"ahler metric ${\mathbb{G}}(\cdot\; ,\cdot)=\Omega({\mathbb{J}}\cdot\; ,\cdot)\;\;$:
\begin{equation}\label{e:metric}
{\mathbb{G}}=-i\left[\frac{1}{(1+\mu_{1}\bar{\mu}_{2})^2}d\mu_{1}\otimes d\bar{\mu}_{2}
                   -\frac{1}{(1+\bar{\mu}_{1}\mu_{2})^2}d\bar{\mu}_{1}\otimes d\mu_{2}\right].
\end{equation}
This metric, which has signature $++--$, is invariant under the action induced on ${\mathbb{L}}({\mathbb{H}}^3)$ by 
the isometry group of ${\mathbb{H}}^3$. Indeed, this has been shown to be the unique K\"ahler metric on 
${\mathbb{L}}({\mathbb{H}}^3)$ with this property \cite{Salvai}.  

In order to transfer geometric data between ${\mathbb{L}}({\mathbb{H}}^3)$ and ${\mathbb{H}}^3$ we use a correspondence 
space:

\begin{figure} [ht]
\vspace*{0.9cm}

\unitlength0.5cm

\begin{picture}(16,6)
\put(5.1,4.8){$\pi_1$}
\put(5.4,6.5){${\mathbb{L}}({\mathbb{H}}^3)\times{\mathbb{R}}$}
\put(7,6){\vector(1,-1){2.8}}
\put(8.5,4.8){$\Phi$}
\put(5.2,2){${\mathbb{L}}({\mathbb{H}}^3)$}
\put(6,6){\vector(0,-1){3}}
\put(10,2.4){${\mathbb{H}}^3$}
\end{picture}
\end{figure}

The key property of this correspondence is that, given $\gamma\in{\mathbb L}({\mathbb H}^3)$, the set 
$\Phi\circ\pi_1^{-1}(\gamma)$ is the oriented geodesic in ${\mathbb H}^3$, while, for a point $p\in{\mathbb H}^3$,
$\pi_1\circ\Phi^{-1}(p)$ is the set of oriented geodesics in ${\mathbb L}({\mathbb H}^3)$ that pass through $p$.

The map $\Phi$ takes an oriented geodesic $\gamma$ in ${\mathbb L}({\mathbb H}^3)$ and a real number $r$ to
the point on  $\gamma$ an affine parameter distance $r$ from some fixed point on the geodesic. This choice of
point on each geodesic can be made globally, but we more often just use a local choice, which is sufficient for our
purposes.

In terms of holomorphic coordinates ($\mu_1,\mu_2$) on ${\mathbb L}({\mathbb H}^3)$ and upper-half space
coordinates ($x^0,x^1,x^2$) the map $\Phi$ has expression:
\[
z=\frac{1-\mu_1\bar{\mu}_2}{2\bar{\mu}_2}+\left(\frac{1+\mu_1\bar{\mu}_2}{2\bar{\mu}_2}\right)\tanh r,
\qquad\qquad 
t=\frac{|1+\bar{\mu}_1\mu_2|}{2|\mu_2|\cosh r}.
\]
where $z=x^1+ix^2$ and $t=x^0$.

We often use local coordinates ($\xi,\eta$) on ${\mathbb L}({\mathbb H}^3)$ which are related to the holomorphic
coordinates by
\begin{equation}\label{e:trans}
\xi=\frac{2\mu_2}{1+\bar{\mu}_1\mu_2},
\qquad\qquad 
\eta=\frac{1-\mu_1\bar{\mu}_2}{2\bar{\mu}_2}.
\end{equation}
The resulting expression for $\Phi$ is
\begin{equation}\label{e:t0z02}
z=\eta+\frac{\tanh r}{\bar{\xi}},
\qquad\qquad 
t=\frac{1}{|\xi|\cosh r}.
\end{equation}

\vspace{0.2in}

\section{Surfaces in ${\mathbb{L}}({\mathbb{H}}^3)$}

In this section we investigate {\it oriented geodesic congruences}, that is, surfaces 
$\Sigma\subset {\mathbb{L}}({\mathbb{H}}^3)$, or two parameter families of oriented geodesics in ${\mathbb{H}}^3$.

\vspace{0.1in}

\subsection{Coordinates and an adapted frame in ${\mathbb{H}}^3$} 

Let $(x^0,x^1,x^2)$ be the standard coordinates on the upper half space model of ${\mathbb{H}}^3$ with $x^0>0$ and set 
$z=x^1+ix^2$, $\bar{z}=x^1-ix^2$ and $t=x^0$.

\begin{Def}
A {\it null frame} in ${\mathbb{H}}^3$ is a trio $\{e_{(0)},e_{(+)},e_{(-)}\}$ of complex vector fields in 
${\mathbb{C}}\otimes $T${\mathbb{H}}^3$ where $e_{(0)}$ is real, $e_{(+)}$ is the complex conjugate of $e_{(-)}$ 
and they satisfy the following properties:
\[
<e_{(0)}\; ,\; e_{(0)}>=1,\quad <e_{(0)}\; ,\; e_{(+)}>=<e_{(+)}\; ,\; e_{(+)}>=0,\quad <e_{(+)}\; ,\; e_{(-)}>=1,
\]
where $<\; ,\; >$ is the hyperbolic metric.
\end{Def}
Given an orthonormal frame $\{e_0,e_1,e_2\}$ on T${\mathbb{H}}^3$ we construct a null frame by
\[
e_{(+)}=\frac{1}{\sqrt{2}}(e_1-ie_2),\qquad e_{(-)}=\frac{1}{\sqrt{2}}(e_1+ie_2),\qquad e_{(0)}=e_0\; .
\]
\begin{Def}
Given a surface $\Sigma\subset {\mathbb{L}}({\mathbb{H}}^3)$ an {\it adapted null frame} is a null frame 
$\{e_{(0)},e_{(+)},e_{(-)}\}$ such that, for each $\gamma\in\Sigma$, we have $e_{(0)}=\dot{\gamma}$, and the orientation 
of $\{e_{(0)},{\mathbb{R}}{\mbox{e}}(e_{(+)}),{\mathbb{I}}{\mbox{m}}(e_{(+)})\}$ is the standard orientation on 
${\mathbb{H}}^3$.
\end{Def}

In what follows we consider a locally parameterized surface $\Sigma\subset {\mathbb{L}}({\mathbb{H}}^3)$ :

\begin{Prop}\label{p:frame}
Let $\Sigma\subset {\mathbb{L}}({\mathbb{H}}^3)$ be a surface given locally by an immersion 
$f:{\mathbb C}\rightarrow {\mathbb{L}}({\mathbb{H}}^3):(\nu,\bar{\nu})\mapsto (\xi(\nu,\bar{\nu}),\eta(\nu,\bar{\nu}))$. Then an adapted null frame $\{e_{(0)},e_{(+)},e_{(-)}\}$ is given by
\[
e_{(0)}=\frac{\partial}{\partial r},\qquad e_{(+)}=a\frac{\partial}{\partial \nu}+b\frac{\partial}{\partial \bar{\nu}}+\omega\frac{\partial}{\partial r},
\]
where
\[
\omega=\sqrt{2}\Delta^{-1}\left(\overline{\partial^{-}F}(\bar{\xi}\;\bar{\partial}\eta+\xi\bar{\partial}\bar{\eta})
        -\overline{\partial^{+}F}(\bar{\xi}\partial\eta+\xi\partial\bar{\eta})\right),
\]
\begin{equation}\label{e:Delta}
a=2\sqrt{2}\Delta^{-1}\;\overline{\partial^{+}F},
\qquad b=-2\sqrt{2}\Delta^{-1}\overline{\partial^{-}F},
\qquad\qquad
\Delta=\partial^{+}F \overline{\partial^{+}F}-\partial^{-}F \overline{\partial^{-}F},
\end{equation}
and
\[
\partial^{+}F=\xi e^{r}\partial\bar{\eta}-\bar{\xi}e^{-r}\partial\eta-e^{r}\partial\ln\xi-e^{-r}\partial\ln\bar{\xi},
\]
\begin{equation}\label{e:d+f}
\partial^{-}F=\xi e^{r}\bar{\partial}\bar{\eta}-\bar{\xi}e^{-r}\bar{\partial}\eta-e^{r}\bar{\partial}\ln\xi
      -e^{-r}\bar{\partial}\ln\bar{\xi}.
\end{equation}
\end{Prop}
\begin{pf}
A parameterized surface $\Sigma$ gives (composing with the map $\Phi$ in (\ref{e:t0z02})) a map 
\begin{equation}\label{e:map}
{\mathbb C}\times {\mathbb{R}}\rightarrow {\mathbb{H}}^3:(\nu,\bar{\nu},r)\mapsto 
     (z(\nu,\bar{\nu},r),\bar{z}(\nu,\bar{\nu},t)).
\end{equation}
Consider now the adapted null frame $\{e_{(0)},e_{(+)},e_{(-)}\}$ given by Proposition 4 of \cite{gag}, 
the derivative of the map (\ref{e:map}) gives 
\[
\frac{\partial}{\partial r}\mapsto e_{(0)},
\]
and
\begin{align}
\frac{\partial}{\partial \nu}&\mapsto {\textstyle{\frac{1}{2}}}\left(\bar{\xi}\partial\eta+\xi\partial \bar{\eta}\right)e_{(0)}
   +{\textstyle{\frac{1}{2\sqrt{2}}}}\left(\xi e^{r}\partial\bar{\eta}-\bar{\xi}e^{-r}\partial\eta
    -e^{r}\partial\ln\xi-e^{-r}\partial\ln\bar{\xi}\right)e_{(+)}\nonumber\\
&\qquad\qquad+{\textstyle{\frac{1}{2\sqrt{2}}}}\left(\bar{\xi} e^{r}\partial\eta-\xi e^{-r}\partial\bar{\eta}
  -e^{r}\partial\ln\bar{\xi}-e^{-r}\partial\ln\xi\right)e_{(-)},\nonumber
\end{align}
The inverse of the above transformation gives the adapted null frame.
\end{pf}

\vspace{0.1in}

\vspace{0.1in}

\subsection{The local geometry of geodesic congruences}

We now investigate the local geometry of surfaces in ${\mathbb L}({\mathbb{H}}^3)$. 
Suppose that $\Sigma$ is such a surface and $\{e_{(0)},e_{(+)},e_{(-)}\}$ is an
adapted null frame, and define the two complex scalars:
\[
\rho=<\nabla_{e_{(-)}}\;e_{(+)},e_{(0)}>,
\qquad\qquad
\sigma=<\nabla_{e_{(+)}}\;e_{(+)},e_{(0)}>.
\]

The complex {\it optical scalar} functions $\rho$ and $\sigma$ describe the first order geometric behaviour of the 
surface. 
The real part of $\rho$ is referred to as the $divergence$, the imaginary part is the $twist$, which is 
denoted by $\lambda$, and $\sigma$ is the $shear$ of the geodesics. They satisfy the Sachs equations
for ${\mathbb{H}}^3$ \cite{PaR}:
\[
\frac{\partial\rho}{\partial r}=\rho^2+\sigma\bar{\sigma}-1,
\qquad\qquad  \frac{\partial\sigma}{\partial r}=(\rho+\bar{\rho})\sigma,
\]
where $r$ is an affine parameter along the geodesics.
Note that if the shear or the twist vanish at some point of the geodesic, they vanish at every point of the geodesic
by iterations of the Sachs equations.

\begin{Prop}
The optical scalars of a surface $(\nu,\bar{\nu})\mapsto (\xi(\nu,\bar{\nu}),\eta(\nu,\bar{\nu}))$
are:

\begin{align}
\sigma&=\frac{2e^{-r}}{\Delta}\left[\left(\xi\bar{\partial}\bar{\eta}+\bar{\partial}\ln \xi\right)\overline{\partial^{-}F}-\left(\xi\partial\bar{\eta}+\partial\ln\xi\right)\overline{\partial^{+}F}\right],\nonumber \\
\rho&=-1+\frac{2e^{-r}}{\Delta}\left[\left(\xi\partial\bar{\eta}+\partial\ln\xi\right)\partial^{-}F-\left(\xi\bar{\partial}\bar{\eta}+\bar{\partial}\ln\xi\right)\partial^{+}F\right],\nonumber
\end{align}
where $\partial^\pm F$ and $\Delta$ are given by (\ref{e:Delta}) and (\ref{e:d+f}), and $\partial$ is differentiation
with respect to $\nu$. 
\end{Prop}

\begin{pf}
The Levi-Civita connection of the hyperbolic metric is torsion-free and so we have the following expression for
the Lie derivative of $e_{(+)}$ in the $e_{(0)}$ direction 
\[
L_{e_{(0)}}e_{(+)}=\bar{\rho}e_{(+)}+\sigma e_{(-)}.
\]
Recall the expression for $e_{(+)}$ given in Proposition \ref{p:frame} and, since $<e_{(0)}\; ,\; e_{(+)}>=0$, we find
\[
\omega=-a<e_{(0)}\;,\;\frac{\partial}{\partial \nu}>-b<e_{(0)}\;,\;\frac{\partial}{\partial \bar{\nu}}>.
\]
Therefore, $e_{(+)}=aZ_{+}+bZ_{-}$ where
\[
Z_{+}=\frac{\partial}{\partial \nu}-<e_{(0)}\;,\;\frac{\partial}{\partial \nu}>\frac{\partial}{\partial r},
\qquad\qquad  Z_{-}=\frac{\partial}{\partial \bar{\nu}}-<e_{(0)}\;,\;\frac{\partial}{\partial \bar{\nu}}>\frac{\partial}{\partial r}\; .
\]
After computing the hyperbolic metric in ($\nu,\bar{\nu},r$) coordinates we find that
\[
<Z_{+}\;,\;Z_{+}>={\textstyle{\frac{1}{4}}}\partial^{+}F\overline{\partial^{-}F},
\qquad\qquad <Z_{+}\;,\;Z_{-}>={\textstyle{\frac{1}{8}}}\left(\partial^{+}F\overline{\partial^{+}F}
                          +\partial^{-}F\overline{\partial^{-}F}\right).
\]
The $e_{(+)}$ component of $L_{e_{(0)}}e_{(+)}$ is $\sigma=<L_{e_{(0)}}e_{(+)}\;,\;e_{(+)}>$
and so
\[
\sigma=\frac{1}{2}\frac{\partial a^2}{\partial r}<Z_{+}\, , \,Z_{+}>+\frac{1}{2}\frac{\partial b^2}{\partial r}<Z_{-}\, , \,Z_{-}>+\frac{\partial ab}{\partial r}<Z_{+}\, , \,Z_{-}>. 
\]
After a lengthy computation we find that 
\[
\sigma=2e^{-r}\frac{\left(\xi\bar{\partial}\bar{\eta}+\bar{\partial}\ln\xi\right)\overline{\partial^{-}F}
      -\left(\xi\partial\bar{\eta}+\partial\ln\xi\right)\overline{\partial^{+}F}}{\partial^{+}F\overline{\partial^{+}F}
                            -\partial^{-}F\overline{\partial^{-}F}}\; .
\]
\vspace{0.1in}
The scalar $\rho$ can be found in a similar way. 
\end{pf}

It's useful to know $\rho$ and $\sigma$ in terms of the holomorphic coordinates $(\mu_1,\mu_2)$. 
\begin{Prop}\label{p:opscal}
Let $\Sigma$ be a geodesic congruence given by an immersion $f:\Sigma\rightarrow {\mathbb{L}}({\mathbb{H}}^3):(\nu,\bar{\nu})\mapsto (\mu_1(\nu,\bar{\nu}),\mu_2(\nu,\bar{\nu}))$. The optical scalars $\rho$ and $\sigma$ are 
\begin{equation}\label{e:sigmarho}
\sigma=\frac{8\mu_2 J_{\bar{2}\;\bar{1}}}{\bar{\mu}_2\Delta |1+\mu_1\bar{\mu}_2|^2},\quad
\rho=-1-\frac{8e^{-r}}{\Delta}\left[\frac{J_{2\bar{1}}}{(1+\bar{\mu}_1\mu_2)^2}e^{r}-\frac{|\mu_2|^2 J_{1\bar{1}}}{|1+\bar{\mu}_1\mu_2|^2}e^{-r}\right],
\end{equation}
where
\[
J_{kl}=\partial\mu_{k}\bar{\partial}\mu_{l}-\bar{\partial}\mu_{k}\partial\mu_{l}\qquad k,l=1,2,\bar{1},\bar{2},
\]
and 
\[
\frac{1}{4}\Delta=\frac{J_{2\bar{2}}}{|\mu_2|^2|1+\mu_1\bar{\mu}_2|^2}e^{2r}+\frac{J_{\bar{2}1}}{(1+\mu_1\bar{\mu}_2)^2}+\frac{J_{\bar{1}2}}{(1+\bar{\mu}_1\mu_2)^2}+\frac{|\mu_2|^2 J_{1\bar{1}}}{|1+\mu_1\bar{\mu}_2|^2}e^{-2r}.
\]
\end{Prop}

\begin{pf}
This follows from the expressions in Proposition \ref{p:opscal} and the transformation (\ref{e:trans}).
\end{pf}

\vspace{0.1in}

\subsection{Lagrangian Surfaces}
\begin{Def}
A point $\gamma$ on a surface $\Sigma\subset{\mathbb{L}}({\mathbb{H}}^3)$ is said to be a {\it Lagrangian point} 
if $f^{*}\Omega=0$ at $\gamma$, where $\Omega$ is the symplectic 2-form on ${\mathbb{L}}({\mathbb{H}}^3)$. A surface
$\Sigma\subset{\mathbb{L}}({\mathbb{H}}^3)$ is said to be  {\it Lagrangian} if all of the points of $\Sigma$ are
Lagrangian points.
\end{Def}

\begin{Prop}\label{p:lagtwist}
A surface $\Sigma$ is Lagrangian iff the imaginary part of $\rho$ (the twist) is zero.
\end{Prop}
\begin{pf}
Let $f\colon\Sigma\rightarrow{\mathbb{L}}({\mathbb{H}}^3)\colon (\nu,\bar{\nu})\mapsto (\mu_1(\nu,\bar{\nu}),\mu_2(\nu,\bar{\nu}))$ be the $C^2$ immersion of the geodesic congruence. 

The twist $\lambda={\mathbb I}\mbox{m}\rho$ is 
\[
\lambda=\frac{4i}{\Delta}\left[\frac{J_{2\bar{1}}}{(1+\bar{\mu}_1\mu_2)^2}-\frac{J_{1\bar{2}}}{(1+\mu_1\bar{\mu}_2)^2}\right].
\]
The symplectic 2-form $\Omega$ in ${\mathbb{L}}({\mathbb{H}}^3)$ is given by (\ref{e:sympl}) and so $f^{*}\Omega$ is 
\[
f^{*}\Omega=\left[\frac{J_{2\bar{1}}}{(1+\bar{\mu}_1\mu_2)^2}-\frac{J_{1\bar{2}}}{(1+\mu_1\bar{\mu}_2)^2}\right]d\nu\wedge d\bar{\nu}=\frac{\Delta}{4i}\lambda\; d\nu\wedge d\bar{\nu}.
\]
The proposition then follows.
\end{pf}
\begin{Def}
An oriented geodesic congruence is $integrable$ iff locally there exists an embedded surface $S$ in ${\mathbb{H}}^3$ such that 
$S$ is orthogonal to the geodesics of the congruence.
\end{Def}
\begin{Prop}
A geodesic congruence $\Sigma$ is integrable iff $\Sigma$ is Lagrangian. 
\end{Prop}
\begin{pf}

The geodesics are integrable iff, in an adapted frame, $\left[e_{(+)},e_{(-)}\right]$ lies in the span of 
$\{e_{(+)},e_{(-)}\}$ (by Frobenius' theorem). But, using the fact that the Levi-Civita connection of the hyperbolic metric
is torsion free,
\[
<\left[e_{(+)},e_{(-)}\right],e_{(0)}>=<\nabla_{e_{(+)}}\;e_{(-)},e_{(0)}>-<\nabla_{e_{(-)}}\;e_{(+)},e_{(0)}>
    =\bar{\rho}-\rho=-2\lambda i
\] 
Thus integrability is equivalent to the imaginary part of $\rho$ vanishing, as claimed.
\end{pf}

\begin{Prop}\label{p:potential}
Let $\Sigma$ be a Lagrangian surface in ${\mathbb L}({\mathbb H}^3)$ parameterized by $\nu\rightarrow(\mu_1(\nu,\bar{\nu}),\mu_2(\nu,\bar{\nu}))$. 
Then the surfaces $S\subset{\mathbb H}^3$ orthogonal to the geodesics of $\Sigma$ are given by the functions $r=r(\nu,\bar{\nu})$ which solve:
\[
2\partial r=\frac{\mu_2}{\bar{\mu}_1\mu_2+1}\left(\partial\bar{\mu}_1+\frac{\partial\mu_2}{\mu_2^2}\right)+\frac{\bar{\mu}_2}{\mu_1\bar{\mu}_2+1}\left(\partial\mu_1+\frac{\partial\bar{\mu}_2}{\bar{\mu}_2^2}\right).
\]
Here the real constant of integration parameterizes the family of parallel surfaces orthogonal to the geodesics.
\end{Prop}
\begin{pf}
A straightforward computation shows that the parametric surface in ${\mathbb{H}}^3$, obtained by inserting $r=r(\nu,\bar{\nu})$ in (\ref{e:t0z02}), 
is orthogonal to the congruence iff the above condition holds.
\end{pf}
\vspace{0.1in}

\begin{Def}
Given an immersion $f:\Sigma\rightarrow{\mathbb L}({\mathbb H}^3)$, consider the map $(\pi\circ f)_*:T\Sigma\rightarrow T{\mathbb P}^1$, where $\pi$ is
projection onto the first factor of ${\mathbb L}({\mathbb H}^3)={\mathbb P}^1\times{\mathbb P}^1-\overline{\Delta}$. The {\it rank} of the
immersion $f$ at a point $\gamma\in\Sigma$ is defined to be the rank of this map at $\gamma$, which can be 0, 1 or 2.
\end{Def}

\vspace{0.1in}

\begin{Note} 
Alternatively, the rank could be defined by projection onto the second factor. However, by reversing the orientation of the geodesics in a geodesic
congruence the two definitions would switch. In what follows we take projection onto the first factor.
 
A rank 0 Lagrangian surface is a horosphere (see Definition \ref{d:horo}).

A rank 2 surface $\Sigma$, parameterized by $\mu_2=\mu_2(\mu_1,\bar{\mu}_1)$, is Lagrangian iff there is a 
smooth real function $r$ such that
\begin{equation}\label{e:rgraph}
2\partial_1 r=\frac{\partial_1\mu_2}{\mu_2(1+\bar{\mu}_1\mu_2)}+\frac{\partial_1\bar{\mu}_2}{\bar{\mu}_2(1+\mu_1\bar{\mu}_2)}
+\frac{\bar{\mu}_2}{1+\mu_1\bar{\mu}_2},
\end{equation}
where here $\partial_1$ is differentiation with respect to $\mu_1$.
\end{Note}

In the Lagrangian case, the functions $\sigma$ and $\rho$ have the following interpretation in terms of the
second fundamental form of the orthogonal surfaces in ${\mathbb H}^3$.

\begin{Prop}\label{p:lag}
Let $S\subset {\mathbb{H}}^3$ be a $C^2$ immersed surface and $\Sigma\subset {\mathbb{L}}({\mathbb{H}}^3)$ 
be the oriented normal geodesics. Then
\begin{equation}\label{e:optical}
|\sigma|=\frac{1}{2}|\lambda_1-\lambda_2|\qquad\qquad \rho=-\frac{1}{2}(\lambda_1+\lambda_2),
\end{equation}
where $\lambda_1$ and $\lambda_2$ are the principal curvatures of $S$.
\end{Prop}
\begin{pf}
Let $S$ be a C$^2$  surface immersed in ${\Bbb{H}}^3$ and let $N_i$ be
the unit normal 1-form. The second fundamental form is a symmetric two tensor on $S$
defined by
\[
H_{ij}=P_i^kP_j^l\nabla_kN_l,
\]
where $\nabla$ is the Levi-Civita connection on ${\Bbb{H}}^3$  and $P_i^j$ is orthogonal projection onto
the tangent space of $S$. Such a symmetric two tensor has two real eigenvalues, $\lambda_1$ and
$\lambda_2$, at each point of $S$. These are called the {\it principal  curvatures} of $S$ and the
associated eigen-directions are called {\it principal directions} for the surface. Let $\{e_{(1)},e_{(2)}\}$ be 
an orthonormal
eigen-basis for $H_{ij}$, and then the null frame is 
$e_{(+)}={\textstyle{\frac{1}{\sqrt{2}}}}(e_{(1)}+ie_{(2)})e^{i\alpha}$ for some angle $\alpha$.
Computing the null frame components of $H_{ij}$
\[
H_{(++)}=<\nabla_{e_{(+)}}\theta^{(0)},e_{(+)}>=\Gamma_{(+)\;\;(+)}^{\;\;\;(0)}=\sigma,
\]
\[
H_{(+-)}=<\nabla_{e_{(+)}}\theta^{(0)},e_{(-)}>=\Gamma_{(+)\;\;(-)}^{\;\;\;(0)}=-\bar{\rho}=-\rho.
\]
In terms of the real basis
\[
\sigma=H_{(++)}={\textstyle{\frac{1}{2}}}\left(H_{(11)}-H_{(22)}\right)e^{2i\alpha} 
   ={\textstyle{\frac{1}{2}}}\left(\lambda_1-\lambda_2\right)e^{2i\alpha},
\]
\[
\rho=-H_{(+-)}=-{\textstyle{\frac{1}{2}}}\left(H_{(11)}+H_{(22)}\right) 
   =-{\textstyle{\frac{1}{2}}}\left(\lambda_1+\lambda_2\right),
\] 
as claimed.
\end{pf}

The Gauss curvature of $S$ is $\kappa=\rho\bar{\rho}-\sigma\bar{\sigma}-1$, and in terms of local coordinates 
($\nu,\bar{\nu}$) this is:
\begin{equation}\label{e:gausscurvature1}
\kappa=\frac{8}{\Delta}\left[\frac{J_{2\bar{1}}}{(1+\bar{\mu}_1\mu_2)^2}
            +\frac{J_{1\bar{2}}}{(1+\mu_1\bar{\mu}_2)^2}\right].
\end{equation}

\vspace{0.1in}

\subsection{Complex Curves}

\begin{Def}
A point $\gamma$ on a surface $\Sigma\subset{\mathbb{L}}({\mathbb{H}}^3)$ is said to be a {\it complex point} 
if the complex structure ${\mathbb{J}}$ acting on ${\mathbb{L}}({\mathbb{H}}^3)$ preserves $T_{\gamma}\Sigma$. A surface
$\Sigma\subset{\mathbb{L}}({\mathbb{H}}^3)$ is said to be a {\it complex curve} if all of the points of $\Sigma$ are
complex points.
\end{Def}
In particular:
\begin{Prop}
A point $\gamma$ on a surface $\Sigma$ is complex iff the shear vanishes along $\gamma$.
\end{Prop}
\begin{pf}
Let $f\colon\Sigma\rightarrow {\mathbb{L}}({\mathbb{H}}^3)\colon (\nu,\bar{\nu})\mapsto (\mu_1(\nu,\bar{\nu}),\mu_2(\nu,\bar{\nu}))$ be the smooth immersion of $\Sigma$ in ${\mathbb{L}}({\mathbb{H}}^3)$. We consider 
the derivative of the immersion 
$df\colon T\Sigma\rightarrow T{\mathbb{L}}({\mathbb{H}}^3)$. Let $j$ be a conformal structure on $\Sigma$ 
compatible with $\nu$.
We define the sections $\delta^{+}f\in \Lambda^{10}(\Sigma)\otimes T^{10}{\mathbb{L}}({\mathbb{H}}^3)$ and 
$\delta^{-}f\in \Lambda^{01}(\Sigma)\otimes T^{10}{\mathbb{L}}({\mathbb{H}}^3)$ by
\[
\delta^{+}f=\frac{1}{2}(df-{\mathbb{J}}\circ df\circ j)\qquad\qquad \delta^{-}f=\frac{1}{2}(df+{\mathbb{J}}\circ df\circ j).
\]
Then $\delta^{+}f\wedge\delta^{-}f\in \Lambda^{2}(\Sigma)\otimes$ det $T^{10}{\mathbb{L}}({\mathbb{H}}^3)$ works out to be
\[
\delta^{+}f\wedge\delta^{-}f=J_{12}d\nu\wedge d\bar{\nu}.
\]
Now, a point is complex iff this 2-form vanishes. On the other hand we have
\[
\sigma=\frac{8\mu_2 J_{\bar{2}\;\bar{1}}}{\bar{\mu}_2\Delta |1+\mu_1\bar{\mu}_2|^2}=\frac{8\mu_2 \overline{J_{12}}}{\bar{\mu}_2\Delta |1+\mu_1\bar{\mu}_2|^2},
\]
and hence the proposition follows.
\end{pf}

\vspace{0.1in}

\subsection{The Induced Metric}

\begin{Thm}
Let $\Sigma$ be a surface in ${\mathbb{L}}({\mathbb{H}}^3)$. The induced metric is Lorentz (degenerate, 
Riemannian) iff $\; |\sigma|^2-\lambda^2>0\; (=0,<0)$, where $\lambda$ and $\sigma$ are the twist and the shear of the 
geodesic congruence $\Sigma$.
\end{Thm}
\begin{pf}
Let $f$ be the smooth immersion like before with local complex coordinate $\nu$. We have seen that the twist and the shear are given by
\[
\sigma=\frac{8\mu_2 J_{\bar{2}\;\bar{1}}}{\bar{\mu}_2\Delta |1+\mu_1\bar{\mu}_2|^2}\qquad\qquad \lambda=\frac{4i}{\Delta}\left[\frac{J_{2\bar{1}}}{(1+\bar{\mu}_1\mu_2)^2}-\frac{J_{1\bar{2}}}{(1+\mu_1\bar{\mu}_2)^2}\right].
\]
Therefore we obtain
\[
|\sigma|^2-\lambda^2=\frac{16}{\Delta^2}\left[\frac{J_{2\bar{1}}^2}{(1+\bar{\mu}_1\mu_2)^4}+\frac{J_{1\bar{2}}^2}{(1+\mu_1\bar{\mu}_2)^4}+\frac{4}{|1+\mu_1\bar{\mu}_2|^4}\left(J_{12}J_{\bar{2}\;\bar{1}}-\frac{1}{2}J_{2\bar{1}}J_{1\bar{2}}\right)\right].
\]
Pulling back the metric ${\mathbb{G}}$, which is given by (\ref{e:metric}), to $\Sigma$ and computing the determinant 
\[
\det [f^{*}{\mathbb{G}}]=-\frac{\Delta^2}{64}(|\sigma|^2-\lambda^2),
\]
and the theorem follows.
\end{pf}

We have the following corollary:

\begin{Cor}
The induced metric on a Lagrangian surface in ${\mathbb{L}}({\mathbb{H}}^3)$ is either Lorentz or degenerate, the 
latter occurring at umbilic points on the orthogonal surfaces in ${\mathbb{H}}^3$.
\end{Cor}
\begin{pf}
By Proposition \ref{p:lagtwist} a Lagrangian surface has $\lambda=0$, and so, by the preceding theorem, $\Sigma$ is 
Lorentz or degenerate, according to whether $\sigma$ is non-zero or zero. By Proposition \ref{p:lag} a point on the
orthogonal surface at which $\sigma=0$ is an umbilic point: $\lambda_1=\lambda_2$. The result follows.
\end{pf}

Similarly, we have:

\begin{Cor}
The induced metric on a complex curve $\Sigma$ in ${\mathbb{L}}({\mathbb{H}}^3)$ is either Riemannian or degenerate, the 
latter occurring at Lagrangian points on $\Sigma$.
\end{Cor}

\vspace{0.2in}

\section{Weingarten surfaces in ${\mathbb{H}}^3$}

\begin{Def}
A surface $S$ in ${\mathbb{H}}^3$ is said to be $Weingarten$ if the eigenvalues $\lambda_1,\lambda_2$ of the 
second fundamental form of $S$ are functionally related \cite{weingarten}:
\[
d\lambda_1\wedge d\lambda_2=0.
\]
\end{Def}
\vspace{0.1in}

we now prove our main result, namely:

\noindent {\bf Main Theorem}:

{\it
Let $S\subset{\mathbb{H}}^3$ be a $C^4$ smooth immersed oriented surface and $\Sigma\subset{\mathbb{L}}({\mathbb{H}}^3)$ 
be the Lagrangian 
surface formed by the oriented geodesics normal to $S$. Assume that the metric ${\mathbb G}_\Sigma$
induced on $\Sigma$ by the neutral K\"ahler metric is non-degenerate. 
Then $S$ is Weingarten iff the Gauss curvature of ${\mathbb G}_\Sigma$ is zero.
}

\begin{pf}

First assume that the normal geodesic congruence $\Sigma$ is of rank 2 and so is locally parameterized by
$\mu_1\mapsto (\mu_1,\mu_2(\mu_1,\bar{\mu}_1))$. 
Our first task is to find the Gauss curvature $K$ of the geodesic congruence. The induced metric $g$, computed by
pulling back (\ref{e:metric}) to $\Sigma$, is:
\begin{equation}\label{e:metric2}
g={\mathbb{G}}_{\Sigma}=2\;{\mathbb{I}}{\mbox{m}}\left[\frac{\partial\bar{\mu}_2 }{(1+\mu_1\bar{\mu}_2)^2}d\mu_1^2
       +\frac{\bar{\partial}\bar{\mu}_2}{(1+\mu_1\bar{\mu}_2)^2}d\mu_1 d\bar{\mu}_1\right],
\end{equation}
where $\partial,\bar{\partial}$ are differentiation with respect of $\mu_1,\bar{\mu}_1$ respectively. 
By Proposition \ref{p:lagtwist} the Lagrangian condition is ${\mathbb{I}}{\mbox{m}}\;\rho=0$ and so by equation
(\ref{e:sigmarho})
\begin{equation}\label{e:lagrangian}
\frac{\bar{\partial}\bar{\mu}_2}{(1+\mu_1\bar{\mu}_2)^2}-\frac{\partial\mu_2}{(1+\bar{\mu}_1\mu_2)^2}=0.
\end{equation}
Thus the induced metric on a Lagrangian surface simplifies to
\[
g=-i\left[\frac{\partial\bar{\mu}_2}{(1+\mu_1\bar{\mu}_2)^2}d\mu_1^2-\frac{\bar{\partial}\mu_2}{(1+\bar{\mu}_1\mu_2)^2}d\bar{\mu}_1^2\right].
\]
The Christoffel symbols $\Gamma_{ij}^{k}$ of the metric ${\mathbb{G}}_{\Sigma}$ are 
\[
\Gamma_{11}^{1}=\frac{1}{2}\partial\left[\ln\left(\frac{\partial\bar{\mu}_2}{(1+\mu_1\bar{\mu}_2)^2}\right)\right],
\qquad \Gamma_{\bar{1}1}^{1}=\frac{1}{2}\bar{\partial}\left[\ln\left(\frac{\partial\bar{\mu}_2}{(1+\mu_1\bar{\mu}_2)^2}\right)\right],
\]
\[
\Gamma_{\bar{1}\;\bar{1}}^{1}=\frac{(1+\mu_1\bar{\mu}_2)^2}{2\partial\bar{\mu}_2}\partial\left[\left(\frac{\bar{\partial}\mu_2}{(1+\bar{\mu}_1\mu_2)^2}\right)\right].
\]
Set
\begin{equation}\label{e:sigma0}
\sigma_0=\frac{\partial\bar{\mu}_2}{(1+\mu_1\bar{\mu}_2)^2},
\end{equation}
so that the induced metric on $\Sigma$ is
\[
g=-i(\sigma_0 d\mu_1^2-\bar{\sigma}_0 d\bar{\mu}_1^{\; 2}),
\] 
and the Christoffel symbols can be written 
\[
\Gamma_{11}^{1}=\frac{\partial\sigma_0}{2\sigma_0},
\qquad\qquad \Gamma_{\bar{1}1}^{1}=\frac{\bar{\partial}\sigma_0}{2\sigma_0},
\qquad\qquad \Gamma_{\bar{1}\;\bar{1}}^{1}=\frac{\partial\bar{\sigma}_0}{2\sigma_0}\; .
\]
The only non vanishing components of the Riemann tensor turn out to be
\[
R_{1\bar{1}1\bar{1}}=\frac{i}{2}(\bar{\partial}^2\sigma_0-\partial^2\bar{\sigma}_0)
+\frac{i}{4}\left[\frac{(\partial\bar{\sigma}_0)^2-\bar{\partial}\bar{\sigma}_0\bar{\partial}\sigma_0}{\bar{\sigma}_0}
-\frac{(\bar{\partial}\sigma_0)^2-\partial\sigma_0\partial\bar{\sigma}_0}{\sigma_0}\right].
\]
For a 2-dimensional surface, we have the following relation that gives the Gauss curvature $K$:
\[
R_{abcd}=K(g_{ac}g_{bd}-g_{ad}g_{bc}).
\]
Using the above relation we obtain
\[
K=\frac{R_{1\bar{1}1\bar{1}}}{g_{11}g_{\bar{1}\;\bar{1}}}=\frac{1}{|\sigma_0|^2}R_{1\bar{1}1\bar{1}}\; ,
\]
and so
\begin{equation}\label{e:curvature1}
K=\frac{i}{4|\sigma_0|^2}\left[2(\bar{\partial}^2\sigma_0-\partial^2\bar{\sigma}_0)+\frac{(\partial\bar{\sigma}_0)^2-\bar{\partial}\bar{\sigma}_0\bar{\partial}\sigma_0}{\bar{\sigma}_0}-\frac{(\bar{\partial}\sigma_0)^2-\partial\sigma_0\partial\bar{\sigma}_0}{\sigma_0}\right].
\end{equation}

If $\sigma_0=0$ then the induced metric vanishes, and we therefore exclude this case and assume that $\sigma_0\neq 0$.

Now set
\begin{equation}\label{e:rho0}
\rho_0=\frac{\partial\mu_2}{(1+\bar{\mu}_1\mu_2)^2}.
\end{equation}
Since $\Sigma$ is Lagrangian, $\rho_0$ is real function. Differentiating (\ref{e:sigma0}) and (\ref{e:rho0}), we obtain 
\[
\bar{\partial}\rho_0=\frac{\partial\bar{\partial}\mu_2}{(1+\bar{\mu}_1\mu_2)^2}-\frac{2\bar{\mu}_1\partial\mu_2\bar{\partial}\mu_2}{(1+\bar{\mu}_1\mu_2)^3}-\frac{2\mu_2\partial\mu_2}{(1+\bar{\mu}_1\mu_2)^3}
\]
\[
\partial\bar{\sigma}_0=\frac{\partial\bar{\partial}\mu_2}{(1+\bar{\mu}_1\mu_2)^2}-\frac{2\bar{\mu}_1\partial\mu_2\bar{\partial}\mu_2}{(1+\bar{\mu}_1\mu_2)^3}, 
\]
and therefore we have established the identity:
\begin{equation}\label{e:codazzi}
\partial\bar{\sigma}_0=\bar{\partial}\rho_0+\frac{2\mu_2}{1+\bar{\mu}_1\mu_2}\rho_0.
\end{equation}
We have seen that the Gauss curvature (\ref{e:curvature1}) has a term with second order derivatives of $\sigma_0$ and $\bar{\sigma}_0$. We can reduce the order of the derivatives using  (\ref{e:codazzi}):
\[
\bar{\partial}^2\sigma_0=\partial\bar{\partial}\rho_0+\bar{\partial}\left(\frac{2\bar{\mu}_2}{1+\mu_1\bar{\mu}_2}\rho_0\right)\qquad\qquad \partial^2\bar{\sigma}_0=\partial\bar{\partial}\rho_0+\partial\left(\frac{2\mu_2}{1+\bar{\mu}_1\mu_2}\rho_0\right),
\]
and subtracting we obtain
\[
\bar{\partial}^2\sigma_0-\partial^2\bar{\sigma}_0=2\left(\frac{\bar{\mu}_2}{1+\mu_1\bar{\mu}_2}\partial\bar{\sigma}_0-\frac{\mu_2}{1+\bar{\mu}_1\mu_2}\bar{\partial}\sigma_0\right).
\]
We now substitute the above into the expression (\ref{e:curvature1}) and obtain the Gauss curvature of $\Sigma$:
\[
K=\frac{i}{|\sigma_0|^2}\left[\frac{\bar{\mu}_2\partial\bar{\sigma}_0}{1+\mu_1\bar{\mu}_2}-\frac{\mu_2\bar{\partial}\sigma_0 }{1+\bar{\mu}_1\mu_2}+\frac{1}{4}\left(\frac{(\partial\bar{\sigma}_0)^2-\bar{\partial}\bar{\sigma}_0\bar{\partial}\sigma_0}{\bar{\sigma}_0}-\frac{(\bar{\partial}\sigma_0)^2-\partial\sigma_0\partial\bar{\sigma}_0}{\sigma_0}\right)\right].
\]

Denote the Gauss curvature of the orthogonal surface $S$ in ${\mathbb{H}}^3$ by $\kappa$. If $\kappa=0$ then we have from equation 
(\ref{e:gausscurvature1}) in graph coordinates that $\partial\mu_2=0$. Then, by expression (\ref{e:rho0}) we have that $\rho_0=0$ and from the identity
(\ref{e:codazzi}) we conclude that $\partial\bar{\sigma}_0=0$. 
Substituting this in the expression above for $K$, we have $K=0$ and the surface in ${\mathbb L}({\mathbb H}^3)$ is scalar flat, as claimed.

Now assume that $\kappa\neq0$. If $\lambda_1,\lambda_2$ are the eigenvalues of the 2nd fundamental form, then $\kappa=\lambda_1\lambda_2-1$. 
From (\ref{e:optical}) we obtain
\[
d\left(\frac{|\sigma|^2}{\kappa^2}\right)\wedge d\left(\frac{\rho+1}{\kappa}\right)=\frac{(\lambda_1-\lambda_2)(\lambda_1-1)(\lambda_2-1)}
     {2(\lambda_1\lambda_2-1)^4}d\lambda_1\wedge d\lambda_2 .
\]
Observe that, on an open set,
\[
d\left(\frac{|\sigma|^2}{\kappa^2}\right)\wedge d\left(\frac{\rho+1}{\kappa}\right)=0\;\; \; \;\; \mbox{iff} \;\;\;\;\; d\lambda_1\wedge d\lambda_2=0,
\]
and therefore a necessary and sufficient condition for the the surface $S$ to be Weingarten is 
\[
d\left(\frac{|\sigma|^2}{\kappa^2}\right)\wedge d\left(\frac{\rho+1}{\kappa}\right)=0.
\]
We are ready now to prove our result.
We need to show that 
\[
d\left(\frac{|\sigma|^2}{\kappa^2}\right)\wedge d\left(\frac{\rho+1}{\kappa}\right)\propto Kd\mu_1\wedge d\bar{\mu}_1 .
\]
Expressing the Gauss curvature $\kappa$ of $S$, and $\rho,\sigma$ in terms of $\rho_0,\sigma_0$
\[
\kappa=\frac{16}{\Delta}\rho_0,
\qquad \rho=-1-\frac{8}{\Delta}\left(\rho_0-\frac{|\mu_2|^2}{|1+\bar{\mu}_1\mu_2|^2}e^{-2r}\right),
\qquad \sigma=\frac{8\mu_2(1+\mu_1\bar{\mu}_2)}{\bar{\mu}_2\Delta(1+\bar{\mu}_1\mu_2)}\sigma_0,
\]
then
\[
\frac{|\sigma|^2}{\kappa^2}=\frac{|\sigma_0|^2}{4\rho_0^2},
\qquad\qquad \frac{\rho+1}{\kappa}=-\frac{1}{2}+\frac{1}{2}\frac{|\mu_2|^2}{|1+\bar{\mu}_1\mu_2|^2}\frac{e^{-2r}}{\rho_0}\; .
\]
Computing the derivatives of the above we obtain
\begin{align}
d\left(\frac{\rho+1}{\kappa}\right)&=-\frac{\rho_0^{-1}|\mu_2|^2 e^{-2r}}{2|1+\bar{\mu}_1\mu_2|^2}\left[\left(\frac{\partial\rho_0}{\rho_0}+\frac{2\bar{\mu}_2}{1+\mu_1\bar{\mu}_2}\right)d\mu_1+\left(\frac{\bar{\partial}\rho_0}{\rho_0}+\frac{2\mu_2}{1+\bar{\mu}_1\mu_2}\right)d\bar{\mu}_1\right],\nonumber \\
d\left(\frac{|\sigma|^2}{\kappa^2}\right)&=
\left[\left(\partial|\sigma_0|^2-2|\sigma_0|^2\frac{\partial\rho_0}{\rho_0}\right)d\mu_1
+\left(\bar{\partial}|\sigma_0|^2-2|\sigma_0|^2\frac{\bar{\partial}\rho_0}{\rho_0}\right)d\bar{\mu}_1\right]
(4\rho_0^2)^{-1}\; ,\nonumber
\end{align}
and hence we find
\[
d\left(\frac{|\sigma|^2}{\kappa^2}\right)\wedge d\left(\frac{\rho+1}{\kappa}\right)=
\frac{|\mu_2|^2|\sigma_0|^4}{2ie^{2r}\rho_0^4|1+\mu_1\bar{\mu}_2|^3}Kd\mu_1\wedge d\bar{\mu}_1\; .
\]
Therefore the surface $S$ is Weingarten iff $K=0$.

Now suppose that  $\Sigma\subset {\mathbb L}({\mathbb H}^3)$ is of rank 1, with immersion 
$f:\Sigma\rightarrow {\mathbb L}({\mathbb H}^3):\nu=u+iv\mapsto (\mu_1(u),\mu_2(u,v))$.
The induced metric $g=f^{\ast}{\mathbb G}$ has components
\begin{align}
g_{uu}&=-i\left[\frac{\partial_{u}\mu_1\partial_{u}\bar{\mu}_2}{(1+\mu_1\bar{\mu}_2)^2}-\frac{\partial_{u}\bar{\mu}_1\partial_{u}\mu_2}{(1+\bar{\mu}_1\mu_2)^2}\right]\nonumber \\
g_{uv}&=g_{vu}=-\frac{i}{2}\left[\frac{\partial_{u}\mu_1\partial_{v}\bar{\mu}_2}{(1+\mu_1\bar{\mu}_2)^2}-\frac{\partial_{u}\bar{\mu}_1\partial_{v}\mu_2}{(1+\bar{\mu}_1\mu_2)^2}\right]\label{e:ajghie} ,
\end{align}
with inverse metric
\begin{equation}\label{e:bkwce}
g^{uu}=0\qquad g^{uv}=g^{vu}=-\frac{1}{g_{uv}}.
\end{equation}

The only non-vanishing Christoffel symbols are:
\begin{align}
\Gamma^{u}_{uu}&=g^{uv}\partial_{u}g_{uv}-\frac{1}{2}g^{uv}\partial_{v}g_{uu}\nonumber \\
\Gamma^{v}_{uu}&=g^{vv}\partial_{u}g_{uv}+\frac{1}{2}g^{uv}\partial_{u}g_{uu}-\frac{1}{2}g^{vv}\partial_{v}g_{uu}\nonumber \\
\Gamma^{v}_{uv}&=\frac{1}{2}g^{uv}\partial_{v}g_{uu}\qquad\qquad \Gamma^{v}_{vv}=g^{uv}\partial_{v}g_{uv} \label{e:gammauuu} 
\end{align}
The Gauss curvature $K$ of the surface $\Sigma$ is given by
\begin{equation}\label{e:koef}
R_{vuvu}=K(g_{uu}g_{vv}-g_{vu}g_{vu})=-g_{vu}^2 K.
\end{equation}
But we know that
\[
R_{vuvu}=g_{vi}R^{i}_{uvu}=g_{vu}R^{u}_{uvu},
\]
and after a brief computation we find that
\begin{equation}\label{e:riemcurv0}
R^{u}_{uvu}=\partial_{v}\Gamma_{uu}^{u}.
\end{equation}

We first find $\Gamma_{uu}^{u}$. From the expressions (\ref{e:gammauuu}) for the Christoffel symbols we need to find 
$\partial_{u}g_{uv}$ and $\partial_{v}g_{uu}$.
\[
\partial_{u}g_{uv}={\mathbb I}{\mbox{m}}\left[\frac{\partial_{u}^2 \mu_1\partial_{v}\bar{\mu}_2}{(1+\mu_1\bar{\mu}_2)^2}+\frac{\partial_{u}\mu_1\partial_{u}\partial_{v}\bar{\mu}_2}{(1+\mu_1\bar{\mu}_2)^2}-\frac{2\mu_1\partial_{u}\mu_1\partial_{u}\bar{\mu}_2\partial_{v}\bar{\mu}_2}{(1+\mu_1\bar{\mu}_2)^3}-\frac{2\bar{\mu}_2(\partial_{u}\mu_1)^2\partial_{v}\bar{\mu}_2}{(1+\mu_1\bar{\mu}_2)^3}\right],
\]
while,
\[
\frac{1}{2}\partial_{v}g_{uu}={\mathbb I}{\mbox{m}}\left[\frac{\partial_{u}\mu_1\partial_{u}\partial_{v}\bar{\mu}_2}{(1+\bar{\mu}_1\mu_2)^2}-\frac{2\mu_1\partial_{u}\mu_1\partial_{u}\bar{\mu}_2\partial_{v}\bar{\mu}_2}{(1+\mu_1\bar{\mu}_2)^3}\right].
\]
Thus 
\[
\partial_{u}g_{uv}-\frac{1}{2}\partial_{v}g_{uu}={\mathbb I}{\mbox{m}}\left[\frac{\partial_{u}^2 \mu_1\partial_{v}\bar{\mu}_2}{(1+\mu_1\bar{\mu}_2)^2}-\frac{2\bar{\mu}_2(\partial_{u}\mu_1)^2\partial_{v}\bar{\mu}_2}{(1+\mu_1\bar{\mu}_2)^3}\right].
\]

Now, the Lagrangian condition is:
\[
\frac{\partial_{u}\mu_1\partial_{v}\bar{\mu}_2}{(1+\mu_1\bar{\mu}_2)^2}=-\frac{\partial_{u}\bar{\mu}_1\partial_{v}\mu_2}{(1+\bar{\mu}_1\mu_2)^2},
\]
which when applied to (\ref{e:bkwce}) gives
\[
g^{uv}=\frac{(1+\mu_1\bar{\mu}_2)^2}{i\partial_{u}\mu_1\partial_{v}\bar{\mu}_2}=-\frac{(1+\bar{\mu}_1\mu_2)^2}{i\partial_{u}\bar{\mu}_1\partial_{v}\mu_2}.
\]
Finally,
\begin{align}
\Gamma_{uu}^{u}&=g^{uv}(\partial_{u}g_{uv}-\frac{1}{2}\partial_{v}g_{uu})\nonumber \\
&=\frac{ig^{uv}}{2}\left[-\frac{\partial_{u}^2 \mu_1\partial_{v}\bar{\mu}_2}{(1+\mu_1\bar{\mu}_2)^2}+\frac{\partial_{u}^2 \bar{\mu}_1\partial_{v}\mu_2}{(1+\bar{\mu}_1\mu_2)^2}+\frac{2\bar{\mu}_2(\partial_{u}\mu_1)^2\partial_{v}\bar{\mu}_2}{(1+\mu_1\bar{\mu}_2)^3}-\frac{2\mu_2(\partial_{u}\bar{\mu}_1)^2\partial_{v}\mu_2}{(1+\bar{\mu}_1\mu_2)^3}\right]\nonumber \\
&={\mathbb I}{\mbox{m}}\left[\frac{\partial_{u}^2 \mu_1\partial_{v}\bar{\mu}_2}{(1+\mu_1\bar{\mu}_2)^2}\frac{(1+\mu_1\bar{\mu}_2)^2}{i\partial_{u}\mu_1\partial_{v}\bar{\mu}_2}+\frac{\partial_{u}^2 \bar{\mu}_1\partial_{v}\mu_2}{(1+\bar{\mu}_1\mu_2)^2}\frac{(1+\bar{\mu}_1\mu_2)^2}{i\partial_{u}\bar{\mu}_1\partial_{v}\mu_2}\right]\nonumber \\
&=-{\mathbb R}{\mbox{e}}\left(\frac{\partial_{u}^2\mu_1}{\partial_{u}\mu_1}-\frac{2\bar{\mu}_2\partial_{u}\mu_1}{1+\mu_1\bar{\mu}_2}\right).\nonumber
\end{align}
Now if we derive the above expression with respect on $v$, we get:
\begin{align}
\partial_{v}\Gamma_{uu}^{u}&=-\frac{1}{2}\partial_{v}\left(\frac{\partial_{u}^2\mu_1}{\partial_{u}\mu_1}+\frac{\partial_{u}^2\bar{\mu}_1}{\partial_{u}\bar{\mu}_1}-\frac{2\bar{\mu}_2\partial_{u}\mu_1}{1+\mu_1\bar{\mu}_2}-\frac{2\mu_2\partial_{u}\bar{\mu}_1}{1+\bar{\mu}_1\mu_2}\right)\nonumber \\
&=\partial_{v}\left(\frac{\bar{\mu}_2\partial_{u}\mu_1}{1+\mu_1\bar{\mu}_2}+\frac{\mu_2\partial_{u}\bar{\mu}_1}{1+\bar{\mu}_1\mu_2}\right)\nonumber \\
&=\frac{\partial_{v}\bar{\mu}_2\partial_{u}\mu_1}{(1+\mu_1\bar{\mu}_2)^2}+\frac{\partial_{v}\mu_2\partial_{u}\bar{\mu}_1}{(1+\bar{\mu}_1\mu_2)^2}\nonumber \\
&=0.\nonumber 
\end{align}
From (\ref{e:riemcurv0}) we have that $R_{uvu}^{u}=\partial_{v}\Gamma_{uu}^{u}=0$, and so $R_{vuvu}=g_{vu}R_{uvu}^{u}=0$. Therefore the Gauss 
curvature $K$ given is
\[
K=-\frac{R_{vuvu}}{g_{uv}^2}=0.
\] 
Thus $\Sigma$ is scalar flat.

We now prove the converse: that if $\Sigma\subset {\mathbb L}({\mathbb H}^3)$ is a Lagrangian surface of rank 1, 
then its orthogonal surface $S\subset {\mathbb H}^3$ is 
Weingarten.

Let $\lambda_1,\lambda_2$ be the principal curvatures of $S$. Then, by Proposition \ref{p:lag}
\begin{equation}\label{e:fwvfejw001}
\kappa=\lambda_1\lambda_2-1\qquad\qquad |\sigma|=\frac{1}{2}|\lambda_1-\lambda_2|,
\end{equation}
where $\kappa$ is the Gauss curvature of $S$ and $\sigma$ is the shear of $\Sigma$.

Since $\Sigma$ is assumed to be of rank 1, it is given by $\mu_1=\mu_1(u)$ and $\mu_2=\mu_2(u,v)$. We know that the shear is
\[
\sigma=\frac{8\mu_2 J_{\bar{2}\;\bar{1}}}{\bar{\mu}_2\Delta |1+\mu_1\bar{\mu}_2|^2},
\]
and therefore,
\begin{equation}\label{e:ijgdfiwf}
|\sigma|^2=\frac{64J_{\bar{2}\;\bar{1}}J_{12}}{\Delta^2 |1+\mu_1\bar{\mu}_2|^4}=\frac{16}{\Delta^2}\frac{\partial_{u}\mu_1\partial_{u}\bar{\mu}_1\partial_{v}\mu_2\partial_{v}\bar{\mu}_2}{|1+\mu_1\bar{\mu}_2|^4}.
\end{equation}
On the other hand, the Gauss curvature $\kappa$ of $S$ is 
\begin{align}
\kappa&=\frac{8}{\Delta}\left[\frac{J_{2\bar{1}}}{(1+\bar{\mu}_1\mu_2)^2}+\frac{J_{1\bar{2}}}{(1+\mu_1\bar{\mu}_2)^2}\right]\nonumber \\
&=\frac{8}{\Delta}\left[-\frac{i\partial_{u}\bar{\mu}_1\partial_{v}\mu_2}{2(1+\bar{\mu}_1\mu_2)^2}+\frac{i\partial_{u}\mu_1\partial_{v}\bar{\mu}_2}{2(1+\mu_1\bar{\mu}_2)^2}\right]\nonumber \\
&=\frac{8i}{\Delta}\frac{\partial_{u}\mu_1\partial_{v}\bar{\mu}_2}{(1+\mu_1\bar{\mu}_2)^2},\nonumber
\end{align}
which  gives
\begin{align}
\kappa^2&=-\frac{64}{\Delta^2}\frac{\partial_{u}\mu_1\partial_{v}\bar{\mu}_2}{(1+\mu_1\bar{\mu}_2)^2}\frac{\partial_{u}\mu_1\partial_{v}\bar{\mu}_2}{(1+\mu_1\bar{\mu}_2)^2}\nonumber \\
&=\frac{64}{\Delta^2}\frac{\partial_{u}\mu_1\partial_{v}\bar{\mu}_2}{(1+\mu_1\bar{\mu}_2)^2}\frac{\partial_{u}\bar{\mu}_1\partial_{v}\mu_2}{(1+\bar{\mu}_1\mu_2)^2}\nonumber \\
&=\frac{64}{\Delta^2}\frac{\partial_{u}\mu_1\partial_{v}\bar{\mu}_2\partial_{u}\bar{\mu}_1\partial_{v}\mu_2}{|1+\mu_1\bar{\mu}_2|^4}\label{e:qsakfo}.
\end{align}

From (\ref{e:ijgdfiwf}) and (\ref{e:qsakfo}) we observe that $\kappa^2=4|\sigma|^2$ and substituting their expressions given in (\ref{e:fwvfejw001}) we obtain
\[
(\lambda_1\lambda_2-1)^2=(\lambda_1-\lambda_2)^2,
\]
or
\[
(\lambda_1-1)(\lambda_1+1)(\lambda_2-1)(\lambda_2+1)=0.
\]
Therefore the surface $S$ is Weingarten and this completes the rank 1 case.

A rank 0 Lagrangian surface is a horosphere, on which the metric is degenerate, and so we have completed the proof of the stated result.

\end{pf}

\vspace{0.2in}

\section{Examples}

\subsection{Totally Null Surfaces}

\begin{Def}
A point $\gamma$ on a surface $\Sigma\subset{\mathbb L}({\mathbb H}^3)$ is a {\it totally null point} if the induced metric on $\Sigma$ is identically
zero at $\gamma$. A surface is {\it totally null} if all of its points are totally null points.
\end{Def}

In \cite{gagak} it is shown that  there are two  types of totally null planes: $\alpha$-planes and $\beta$-planes. On the former, the anti-self-dual 
2-forms vanish, while on the latter the self-dual 2-forms vanish. In particular, the $\alpha$-planes are holomorphic and Lagrangian. 

We now classify the totally null surfaces in ${\mathbb L}({\mathbb H}^3)$. Before doing so, let us introduce some 
terminology.

\begin{Def}
Given a point $p\in{\mathbb{H}}^3$, the {\it geodesic sphere of radius $r_0$ and centre p} is 
the surface in ${\mathbb{H}}^3$ obtained by flowing a distance $r_0$ along all oriented geodesics through p. 
\end{Def}

The geodesic spheres have the following useful characterization:

\begin{Prop}
The oriented normals to a geodesic sphere with centre $(t_0,z_0)$ in the upper half space model 
are given by the solution set of:
\begin{equation}\label{e:sphere}
\bar{z}_0\mu_1\mu_2+(t_0^2+z_0\bar{z}_0)\mu_2-\mu_1-z_0=0.
\end{equation}
\end{Prop}
\begin{pf}
Let $(\xi,\eta)$ be an oriented geodesic which is normal to a geodesic sphere in ${\mathbb{H}}^3$ 
with centre $(t_0,z_0)$. Then there exists $r_0\in {\mathbb{R}}$ such that (see equation (\ref{e:t0z02}))
\begin{equation}\label{e:t0z01}
t_0=\frac{|1+\bar{\mu}_1\mu_2|}{2|\mu_2|\cosh r_0}\qquad\qquad z_0=\frac{1-\mu_1\bar{\mu}_2}{2\bar{\mu}_2}+\left(\frac{1+\mu_1\bar{\mu}_2}{2\bar{\mu}_2}\right)\tanh r_0\; .
\end{equation}
It is now easy to show that (\ref{e:t0z01}) implies (\ref{e:sphere})

Conversely, suppose (\ref{e:sphere}) holds for an oriented geodesic $\gamma=(\mu_1,\mu_2)$. Using the coordinate change 
(\ref{e:trans}) from $(\mu_1,\mu_2)$ to $(\xi,\eta)$ the relation (\ref{e:sphere}) becomes:
\[
\xi(\bar{z}_0-\bar{\eta})-\bar{\xi}(z_0-\eta)+|\xi|^2\left(t_0^2+|z_0-\eta|^2-\frac{1}{|\xi|^2}\right)=0\; .
\]
Since
\[
\xi(\bar{z}_0-\bar{\eta})-\bar{\xi}(z_0-\eta)\in i{\mathbb{R}}\qquad{\mbox{ and }}\qquad |\xi|^2\left(t_0^2+|z_0-\eta|^2-\frac{1}{|\xi|^2}\right)\in {\mathbb{R}},
\]
we obtain
\begin{equation}\label{e:sxesi1}
\xi(\bar{z}_0-\bar{\eta})=\bar{\xi}(z_0-\eta) 
\end{equation}
\begin{equation}\label{e:sxesi2}
t_0^2+|z_0-\eta|^2-\frac{1}{|\xi|^2}=0.
\end{equation}
Thus, by equation (\ref{e:sxesi1}), $\bar{\xi}(z_0-\eta)$ is real, and, by the relation (\ref{e:sxesi2}) (using the fact $t_0>0$), we have
\[
[\bar{\xi}(z_0-\eta)]^2<1,
\]
and therefore there is $r_0\in {\mathbb{R}}$ such that 
\begin{equation}\label{e:sxesi3}
\bar{\xi}(z_0-\eta)=\tanh r_0.
\end{equation}
Substituting (\ref{e:sxesi3}) into (\ref{e:sxesi2}), and using again the fact $t_0>0$ we obtain
\[
t_0=\frac{1}{|\xi|\cosh r_0}\quad{\mbox{ or }}\quad z_0=\eta+\frac{\tanh r_0}{\bar{\xi}}.
\]
Therefore $(t_0,z_0)\in\gamma$ (see equation (\ref{e:t0z02})) and hence  $\gamma$ belongs to the sphere in ${\mathbb{L}}({\mathbb{H}}^3)$ with centre $(t_0,z_0)$. 
\end{pf}

If we fix a point $\mu_0$ on ${\mathbb P}^1$ (considered as the boundary of the ball model of ${\mathbb H}^3$) 
and look at all of the oriented geodesics that end at $\mu_0$, we obtain a surface in 
${\mathbb L}({\mathbb H}^3)$ that is readily found to be Lagrangian. The orthogonal surfaces in ${\mathbb H}^3$
are well-known:

\begin{Def}\label{d:horo}
A {\it horosphere} is a surface in ${\mathbb H}^3$ whose oriented normals end (or begin) at the same point on
${\mathbb P}^1$.
\end{Def}

\vspace{0.1in}

\subsection{$\alpha$-surfaces}

We now establish:

\vspace{0.1in}

\begin{Thm}
An immersed surface $\Sigma$ in ${\mathbb L}({\mathbb H}^3)$ is an $\alpha$-surface iff $\Sigma$ is the oriented normal congruence of
\begin{enumerate}
\item a geodesic sphere, or 
\item a horosphere, or
\item a totally geodesic surface
\end{enumerate}
in ${\mathbb{H}}^3$.
\end{Thm}
\begin{pf}
Straight-forward calculations verify that the oriented normal congruence to the three listed classes of surfaces are 
holomorphic and Lagrangian.

Conversely, assume that $\Sigma$ is a holomorphic and Lagrangian. If it is of rank 0, then it is a horosphere. 

Suppose then that the surface is of rank 1.
Thus locally the surface can be parameterized by $\nu=u+iv\rightarrow(\mu_1(u),\mu_2(u,v))$.
A straightforward computation shows that 
\[
J_{\bar{2}\bar{1}}=-{\textstyle{\frac{i}{2}}}\partial_u\bar{\mu}_1\partial_v\bar{\mu}_2
\]
Holomorphicity implies that $J_{\bar{2}\bar{1}}=0$, and so either $\partial_u\bar{\mu}_1=0$ or $\partial_v\bar{\mu}_2=0$. Thus, either the surface
is not of rank 1, or it is not immersed, respectively.

Now assume that the surface is of rank 2 and parameterize the surface $\Sigma$ with
$\mu_1$. The vanishing of the shear implies that $\mu_2$ is a holomorphic function of $\mu_1$. 
Thus, about any given $\gamma\in\Sigma$, it can be expanded in a power series:
\begin{equation}\label{e:series}
\mu_2=\sum_{n=0}^{\infty}A_{n}\mu_1^{n},
\end{equation}
where $A_{n}$ are complex numbers. The Lagrangian condition (\ref{e:lagrangian}) says that
\[
(1+\mu_1\bar{\mu}_2)^2\partial\mu_2=(1+\bar{\mu}_1\mu_2)^2\bar{\partial}\bar{\mu_2}.
\]
Inserting the series (\ref{e:series}) in this, the holomorphic terms lead to the following relation:
\[
A_1-\bar{A}_1+\sum_{n=1}^{\infty}\left[(n+1)A_{n+1}+2nA_{n}\bar{A}_0+(n-1)A_{n-1}\bar{A}_0^2\right]\mu_1^{n}=0.
\]
Thus $A_1$ is real and we get the recursion relations
\[
(n+1)A_{n+1}+2nA_{n}\bar{A}_0+(n-1)A_{n-1}\bar{A}_0^2=0\quad n=1,2,3,...
\]
It can be easily proven by induction that this is equivalent to
\[
A_{n}=(-1)^{n-1}\bar{A}_0^{n-1}A_1,
\] 
for any $n=1,2,3,...$, and therefore
\begin{equation}\label{e:olom}
\mu_2=A_0+\sum_{n=1}^{\infty}(-1)^{n-1}\bar{A}_0^{n-1}A_1\mu_1^{n}=\frac{A_0+(A_0\bar{A}_0+A_1)\mu_1}{1+\bar{A}_0\mu_1},
\end{equation}
where $A_0\in {\mathbb{C}}$ and $A_1\in {\mathbb{R}}$.

If $A_1>0$, $\Sigma$ is a geodesic sphere with centre $(z_0,t_0)$ with
\[
z_0=\frac{A_0}{A_1+A_0\bar{A}_0}
\qquad\qquad 
t_0=\frac{\sqrt{A_1}}{A_1+A_0\bar{A}_0},
\]
as can be seen by inserting (\ref{e:olom}) in (\ref{e:sphere}).

On the other hand, if $A_1=0$, $\Sigma$ is clearly a horosphere because $\mu_2=A_0$. 

We now prove that if $A_1<0$ then there exists an orthogonal surface in ${\mathbb{H}}^3$ which is totally geodesic.

For the Lagrangian geodesic congruence given by (\ref{e:olom}), the orthogonal surfaces in ${\mathbb{H}}^3$
are obtained by integrating (\ref{e:rgraph}), which yields
\begin{equation}\label{e:r1}
2r=\ln|A_0+(A_1+A_0\bar{A}_0)\mu_1|^2+C,
\end{equation}
where $C\in {\mathbb{R}}$. The divergence can now be computed by (\ref{e:sigmarho}) and the result is:
\[
\rho=-\frac{A_1 e^{C}+1}{A_1 e^{C}-1}.
\]
Thus the orthogonal surface in ${\mathbb{H}}^3$ obtained by setting $C=-\ln(-A_1)$ in (\ref{e:r1}) is totally geodesic, 
since both $\sigma$ and $\rho$ vanish.
\end{pf}

{\noindent}{\bf Note}:
In proving the preceding, we have established the well known classification of totally umbilic surfaces in ${\mathbb H}^3$
\cite{spivak}.

\vspace{0.1in}

\subsection{$\beta$-surfaces}

In \cite{gagak} the $\beta$-surfaces are classified:

\begin{Thm}
A $\beta$-surface in ${\mathbb{L}}({\mathbb{H}}^3)$ is a piece of a torus which, up to 
isometry, is either
\begin{enumerate}
\item ${\mathbb{L}}({\mathbb{H}}^2)$, where ${\mathbb{H}}^2\subset{\mathbb{H}}^3$, or
\item ${\cal C}_1\times {{\cal C}}_2\subset\; S^2\times S^2-\bar{\Delta} $, where ${\cal C}_1$ is a circle given by the
intersection of the 2-sphere and a plane containing the north pole, and ${{\cal C}}_2$ is the image of ${\cal C}_1$ under 
reflection in the horizontal plane through the origin.
\end{enumerate}
\end{Thm}

\vspace{0.1in}

\subsection{The Holomorphic Structure of ${\mathbb{L}}({\mathbb H}^3)$}\label{s:hol}

Consider ${\mathbb{P}}^1\times {\mathbb{P}}^1$. Then the homology group
$H_{2}({\mathbb{P}}^1\times {\mathbb{P}}^1)\cong {\mathbb{Z}}\oplus {\mathbb{Z}}$ has generators $[h]$ and $[v]$,
where $h$ and $v$ are horizontal and vertical fibres of the natural projections. 
The intersection pairing on these generators is obviously
\[
[h]\cdot \mbox{[v]}=1\qquad\qquad [h]\cdot [h]=\mbox{[v]}\cdot \mbox{[v]}=0.
\]
It is also not hard to show that $[\Delta]=[h]+\mbox{[v]}$ and $[\overline{\Delta}]=[h]-\mbox{[v]}$. 

Let $C$ be a closed holomorphic curve in ${\mathbb{P}}^1\times {\mathbb{P}}^1$. Then
\[
[C]=m[h]+n\mbox{[v]},
\]
where $m\geq 0,\; n\geq 0$ and, with $m$ and $n$ not both zero. Clearly 
\[
[C]\cdot [\Delta]=m+n>0\; ,
\]
so every closed holomorphic curve in ${\mathbb{P}}^1\times {\mathbb{P}}^1$ must intersect the diagonal. We conclude that 
there are no closed holomorphic curves in ${\mathbb{P}}^1\times {\mathbb{P}}^1-\Delta$. If ${\mathbb{L}}({\mathbb{H}}^3)$ 
were biholomorphic to ${\mathbb{P}}^1\times {\mathbb{P}}^1-\Delta$, then by Proposition \ref{p:lag}, there would be no 
closed 
totally umbilic surfaces in ${\mathbb{H}}^3$, which we have just shown is not true (the geodesic spheres 
being the counter-examples). No such restriction applies 
to ${\mathbb{P}}^1\times {\mathbb{P}}^1-\overline{\Delta}$, since there are many closed holomorphic curves in 
${\mathbb{P}}^1\times {\mathbb{P}}^1$ that do not intersect $\overline{\Delta}$, one example being $\Delta$:
\[
[\Delta]\cdot [\overline{\Delta}]=1-1=0.
\]
Thus, the distinction between ${\mathbb{P}}^1\times {\mathbb{P}}^1-\Delta$ and 
${\mathbb{P}}^1\times {\mathbb{P}}^1-\overline{\Delta}$ is essential.

\vspace{0.1in}

\subsection{Flat Surfaces}

\begin{Def}
A surface $S$ in ${\mathbb{H}}^3$ is {\it flat} if the Gauss curvature $\kappa$ of $S$ is zero.
\end{Def}

\begin{Thm}
Let $S$ be an oriented $C^2$ smooth immersed surface in ${\mathbb{H}}^3$ with normal geodesic congruence 
$f:\Sigma\rightarrow {\mathbb{L}}({\mathbb{H}}^3)$ parameterized by $\mu_2=\mu_2(\mu_1,\bar{\mu}_1)$. 

Then S is flat iff $\mu_2$ is an anti-holomorphic function of $\mu_1$.
\end{Thm}

\begin{pf}
In this case, the Lagrangian condition for $\Sigma$ is
\[
\frac{\partial\mu_2}{(1+\bar{\mu}_1\mu_2)^2}=\frac{\bar{\partial}\bar{\mu}_2}{(1+\mu_1\bar{\mu}_2)^2},
\]
and then, by equation (\ref{e:gausscurvature1}), the Gauss curvature $\kappa$ of $S$ 
\[
\kappa=\frac{16}{\Delta}\frac{\partial\mu_2}{(1+\bar{\mu}_1\mu_2)^2}.
\]
Therefore the vanishing of the Gauss curvature is equivalent to $\mu_2$ being an anti-holomorphic function. We can see 
as well that if a Lagrangian geodesic congruence has an orthogonal flat surface in ${\mathbb{H}}^3$ then all of
its orthogonal surfaces are flat.

\end{pf}

\vspace{0.1in}

\subsection{Surfaces of Constant Mean Curvature 1}

\begin{Def}
A surface is $S$ in ${\mathbb{H}}^3$ is of {\it constant mean curvature} 1 (CMC 1) if the mean curvature of $S$ is 
equal to 1. Equivalently, the divergence is $\rho=-1$.
\end{Def}

We now prove 
\begin{Thm}
An immersed Lagrangian surface given by $\mu_2=\mu_2(\mu_1,\bar{\mu}_1)$ has an orthogonal CMC 1 surface if and only if
$\sigma_0$ is holomorphic, where
\[
\sigma_0=\frac{\partial\bar{\mu}_2}{(1+\mu_1\bar{\mu}_2)^2}.
\]
\end{Thm}
\begin{pf}
Let $S$ be an oriented $C^2$ smooth surface with constant mean curvature $H$=1. Let 
$\Sigma\subset{\mathbb{L}}({\mathbb{H}}^3)$ be the oriented normal geodesic congruence to $S$. 

Since $\Sigma$ is Lagrangian,
by Proposition \ref{p:potential} 
there exists a function $r:\Sigma\rightarrow {\mathbb{R}}$ satisfying (\ref{e:rgraph}). 

Now, being CMC 1 implies, from $\rho=-1$ and the expression (\ref{e:sigmarho}),
\[
\frac{\partial\mu_2}{(1+\bar{\mu}_1\mu_2)^2}-\frac{|\mu_2|^2}{|1+\bar{\mu}_1\mu_2|^2}e^{-2r}=0,
\]
which can be integrated to
\[
r=-\frac{1}{2}\ln\left(\frac{|1+\bar{\mu}_1\mu_2|^2}{|\mu_2|^2}\rho_0\right),
\]
where $\rho_0$ is defined in equation (\ref{e:rho0}).

Inserting this in (\ref{e:rgraph}) we obtain
\begin{equation}\label{e:conditioncmc}
-\partial\ln\left(\frac{|1+\bar{\mu}_1\mu_2|^2}{|\mu_2|^2}\rho_0\right) =\frac{\partial\mu_2}{\mu_2(1+\bar{\mu}_1\mu_2)}+\frac{\partial\bar{\mu}_2}{\bar{\mu}_2(1+\mu_1\bar{\mu}_2)}+\frac{\bar{\mu}_2}{1+\mu_1\bar{\mu}_2}.
\end{equation}

On the other hand, if a given Lagrangian surface $\Sigma$ satisfies the condition (\ref{e:conditioncmc}), then 
from (\ref{e:rgraph}) we have
\[
2\partial r=-\partial\ln\left(\frac{|1+\bar{\mu}_1\mu_2|^2}{|\mu_2|^2}\rho_0\right),
\]
and hence
\[
r=-\frac{1}{2}\ln\left(\frac{|1+\bar{\mu}_1\mu_2|^2}{|\mu_2|^2}\rho_0\right)+c,
\]
$c$ is real constant. We thus obtain an orthogonal surface $S$ to the geodesic congruence $\Sigma$, parameterized by $r$ with $c=0$. This surface is CMC 1. 

Therefore a Lagrangian geodesic congruence $\mu_2=\mu_2(\mu_1,\bar{\mu}_1)$ has an orthogonal CMC 1 surface if and only if (\ref{e:conditioncmc}) is satisfied. Now we prove that (\ref{e:conditioncmc}) is equivalent to $\bar{\partial}\sigma_0=0$.

The relation (\ref{e:conditioncmc}) can be written
\[
\frac{1+\bar{\mu}_1\mu_2}{\mu_2}\rho_0+\frac{1+\mu_1\bar{\mu}_2}{\bar{\mu_2}}\bar{\sigma}_0+\frac{\bar{\mu}_2}{1+\mu_1\bar{\mu}_2}=-\frac{\bar{\mu}_1\partial\mu_2}{1+\bar{\mu}_1\mu_2}-\frac{\bar{\mu}_2+\mu_1\partial\bar{\mu}_2}{1+\mu_1\bar{\mu}_2}-\frac{\partial\rho_0}{\rho_0}+\frac{\partial\mu_2}{\mu_2}+\frac{\partial\bar{\mu}_2}{\bar{\mu}_2}.
\]
Using (\ref{e:codazzi}) we get
\[
\frac{1+\bar{\mu}_1\mu_2}{\mu_2}\rho_0+\frac{1+\mu_1\bar{\mu}_2}{\bar{\mu_2}}\bar{\sigma}_0+\frac{\bar{\mu}_2}{1+\mu_1\bar{\mu}_2}=-\bar{\mu}_1(1+\bar{\mu}_1\mu_2)\rho_0-\mu_1(1+\mu_1\bar{\mu}_2)\bar{\sigma}_0
\]
\[
-\frac{\bar{\mu}_2}{1+\mu_1\bar{\mu}_2}-\frac{1}{\rho_0}\left(\bar{\partial}\sigma_0-\frac{2\bar{\mu}_2}{1+\mu_1\bar{\mu}_2}\rho_0\right)+\frac{\partial\mu_2}{\mu_2}+\frac{\partial\bar{\mu}_2}{\bar{\mu}_2}.
\]
Cancelling terms using the definitions (\ref{e:sigma0}) and (\ref{e:rho0}),  we find that $\bar{\partial}\sigma_0=0$.
\end{pf}

\vspace{0.2in}
\noindent{\bf Acknowledgement}:

The authors would like to thank Madeeha Khalid, Wilhelm Klingenberg and Jos\'e Ram\'on Mar\'i for many stimulating and
helpful conversations. This work was supported by the Research in Pairs Programme of the 
Mathematisches Forschungsinstitut Oberwolfach, Germany and the Science Foundation Ireland Research Frontiers Programme.

\end{document}